\documentclass[12pt,reqno]{amsart}
\usepackage{amssymb,amsmath,amsthm,fullpage}
\usepackage{pdfsync}

\newtheorem{theorem}{Theorem}
\newtheorem{lemma}[theorem]{Lemma}
\newtheorem{prop}[theorem]{Proposition}

\newtheorem{corollary}[theorem]{Corollary}

\theoremstyle{definition}
\newtheorem{definition}[theorem]{Definition}

\def\eps{\epsilon}

\def\th{\theta}

\DeclareMathOperator{\Dom}{Dom}
\DeclareMathOperator{\Ker}{Ker}
\DeclareMathOperator{\Rre}{Re}
\DeclareMathOperator{\Imm}{Im}
\def\bar{\overline}
\def\w{\widehat}
\def\n{\noindent}
\def\m{\medskip}
\def\eps{\epsilon}

\newcommand{\R}{\mathbb{R}}
\newcommand{\C}{\mathbb{C}}
\newcommand{\Z}{\mathbb{Z}}

\newcommand{\mulj}{\mu^\lambda_j}

\newcommand{\muj}{\mu_j}

\newcommand{\I}{\mathcal{I}}

\newcommand{\z}{\bar z}

\begin{document}

\title{Fundamental Solutions to $\Box_b$ on Certain Quadrics}

\begin{abstract} The purpose of this article is to expand the number of examples for which the complex Green operator, that is, the fundamental solution to the Kohn Laplacian, can
be computed. We use the Lie group structure of quadric submanifolds of $\C^n\times\C^m$ and the group Fourier transform to reduce the $\Box_b$ equation to ones
that can be solved using modified Hermite functions. We use Mehler's formula and investigate  1) quadric hypersurfaces,
where the eigenvalues of the Levi form are not identical
(including possibly zero eigenvalues), and 2) the canonical quadrics
in $\C^4$ of codimension two.
\end{abstract}

\author{Albert Boggess and Andrew Raich}

\thanks{The second author is partially supported by NSF grant DMS-0855822.}

\address{
Department of Mathematics\\ Texas A\&M University\\ Mailstop 3368  \\ College Station, TX
77845-3368 }
\address{
Department of Mathematical Sciences \\ 1 University of Arkansas \\ SCEN 327 \\ Fayetteville, AR 72701}


\email{boggess@math.tamu.edu, araich@uark.edu}

\date{\today}

\maketitle

%
%
\section{Introduction}

The Heisenberg group in $\C^n$ is the primary example
for which the fundamental solution to the Kohn Laplacian $\Box_b$ can be explicitly 
written down. The purpose of this article is to expand 
this library of examples to include: 1) quadric hypersurfaces,
where the eigenvalues of the Levi form are not identical
(including possibly zero eigenvalues), and 2) the canonical quadrics
in $\C^4$ of codimension two. It should be noted that $\Box_b$
does not transform well under biholomorphic changes of coordinates.
In particular, the fundamental solution to $\Box_b$ in the case of 
a strictly pseudoconvex quadric hypersurface where the 
eigenvalues are different cannot be obtained simply 
by rescaling the variables in the Heisenberg case. 

The techniques used in this paper involve representation theory
to reduce $\Box_b$ to a modified Hermite equation, where 
a spectral expansion by Hermite functions provides a quick
answer to the solvability of $\Box_b$. This approach has
been used by many authors (see for example, Peloso and Ricci \cite{PeRi03} and the references therein). Our work then goes one
step further to compute a (nearly) closed form expression 
for the fundamental solution to $\Box_b$ for the example
quadrics, mentioned above, by using this spectral expansion
together with a modification of Mehler's formula. The formulas
obtained clearly exhibit the solution operator for $\Box_b$ in terms
of the eigenvalues of the Levi form and should make it easier
to obtain more precise estimates in terms of these eigenvalues. 

In \cite{FoSt74e}, Folland and Stein give a formula for the fundamental solution of a family of second order operators $\mathcal L_\alpha$ that include the fundamental solution
to the Kohn Laplacian at every form level (for which a solution exists). Since Folland and Stein simply present their solution and demonstrate that it works, it is 
not adaptable to many circumstances. On the other hand, in \cite{PeRi03} Peloso and Ricci give a systematic treatment of $\Box_b$ on quadrics using representation theory and Hermite
expansions. Our motivation is to use the Peloso/Ricci approach but generate formulas in the spirit of \cite{FoSt74e}. 
We realized from our earlier work on the Heisenberg group and quadric submanifolds
\cite{BoRa09, BoRa11} that we could adopt Mehler's formula, one of the workhorse equations with Hermite functions and the $\Box_b$-heat equation, to the $\Box_b$ equation.

Previously to our work, authors attempted to find formulas for the fundamental solution to  $\Box_b$ via the $\Box_b$-heat equation or Hamilton-Jacobi equations. 
In particular, since $e^{-s\Box_b}$ solves
the $\Box_b$-heat equation, if $P$ is the projection onto $\Ker(\Box_b)$, then
\[
\int_0^\infty e^{-s\Box_b}(I-P)\, ds
\]
is the fundamental solution to $\Box_b$. The difficulty with this approach is that the techniques to solve for $e^{-s\Box_b}$ do so up to a partial Fourier transform in the
totally real tangent direction. Calin et.\ al.\  are able to recover the solution for the Heisenberg group and find some formulas for a more general class of quadrics \cite{CaChTi06}. The authors
are able to use this technique to prove estimates on the $\Box_b$-heat kernel on certain decoupled polynomial models in $\C^n$ \cite{Rai06h,Rai07,Rai12h,BoRa11}.
Beals et.\ al.\ use an approach through Hamilton-Jacobi equations to generate integral formulas for fundamental solutions for certain subelliptic equations \cite{BeGaGr96}. Earlier
work by Beals and Greiner in this direction used pseudodifferential operators \cite{BeGr88}.
  
\section{Statement of Results} 
In this section, we summarize the type of results we can obtain from
our techniques. A precise definition of $\Box_b$ and further background
is given in the following sections.

\subsection{Hypersurface Results.}
We consider the case of 
a quadric hypersurface in $\C^{n+1}$ of the form
\[
M= \{(z,w) \in \C^n \times \C; \ \Imm { w}=\sum_{j=1}^n \sigma_j |z_j|^2\},
\ \ \sigma_j \in \R, \ 1 \leq j \leq n.
\]
We shall identify a point $(z,w) \in M$ with the point $(z,t= \Rre w)$
lying in the tangent space. 
The case where all the $\sigma_j$ are equal and nonzero is the 
Heisenberg group, and the fundamental solution to $\Box_b$ 
on $(0,q)$ forms for $1 \leq q \leq n-1$ is well known in this case
(\cite{FoSt74p,FoSt74e,Ste93}). To simplify the notation of our result, we will (mostly) restrict ourselves
to $n=2$.

\begin{theorem}
\label{mixedeigen}
Suppose $n=2$, and $M$ is the hypersurface presented as above
with 
$\sigma_1, \ \sigma_2 >0$. Then, on the space of differential $(0,1)$ forms,
spanned by $d \bar z_1$, $\Box_b$ is solvable 
via a convolution
with a kernel
$N(z,t)$ where
\begin{align*}
\label{hyperformula}
N(z,t)&=\frac{1}{\pi^3} \int_0^1 \frac{\sigma_1 \sigma_2}{(it+s_1(r)\sigma_1|z_1|^2+
s_2(r)\sigma_2|z_2|^2)^2}
\frac{r^{\sigma_1 -1}\, dr}{(1-r^{\sigma_1})(1-r^{\sigma_2})} \\
&+ \frac{1}{\pi^3} \int_0^1 \frac{\sigma_1\sigma_2}{(-it+s_1(r)\sigma_1|z_1|^2
+s_2(r)\sigma_2|z_2|^2)^2}
\frac{r^{\sigma_2-1}\, dr}{(1-r^{\sigma_1})(1-r^{\sigma_2})} 
\end{align*}
where 
\[
s_j(r) = \frac{1+r^{\sigma_j}}{1-r^{\sigma_j}} \ \ \textrm{for} \ j=1, \ 2.
\]
In the case of forms spanned by $d \bar z_2$, the factor of $r^{\sigma_1-1}$
in the numerator is replaced by $r^{\sigma_2-1}$, and vice versa.
\end{theorem}

The convolution mentioned in the above theorem is taken with respect
to the group structure on the quadric hypersurface, which is discussed
in more detail in the following section.
In the case where $\sigma_1=\sigma_2$ an easy conformal rescaling (which does
preserve $\Box_b$) can be used to arrange $\sigma_1=\sigma_2 =1$. In this case,
the change of variables 
\[
s=\frac{1+r}{1-r}, \ \ ds = \frac{2dr}{(1-r)^2}
\]
allows one to compute this integral explicitly to obtain 
\[
N(z,t)= \frac{1}{\pi^3 (|z|^2+it)(|z|^2-it)} = \frac{1}{\pi^3 (|z|^4 +|t|^2)},
\]
the well-known fundamental solution for the Heisenberg group in $\C^3$
on $(0,1)$-forms.

Although we have specialized to the case where $n=2$ and $\sigma_j>0$,
it is clear from the presentation below that the same techniques can handle
$n \geq 2$ and $\sigma_j \in \R$. If one or more of the $\sigma_j$ vanish, 
then the fundamental solution to $\Box_b$ takes a somewhat different form.
To simplify the notation, we consider the case $n=3$ with one vanishing eigenvalue.

\begin{theorem}
\label{zeroeigen}
Suppose $n=3$, and $M$ is the hypersurface presented as above
with 
$\sigma_1= \sigma_2 =1$ and $\sigma_3=0$. Then, on the space of differential $(0,1)$ forms,
spanned by $d \bar z_1$ and $d \bar z_2$, $\Box_b$ is solvable 
via a convolution
with a kernel
$N(z,t)$ where
\[
N(z,t) = \frac{8}{\pi^4}
\int_0^1 \frac{dr}{ |\ln r| (1-r)^2} \Rre \left\{
\frac{1}{\left[ |z'|^2 \left(\frac{1+r}{1-r}\right) + |z_3|^2 \frac{2}{ |\ln r|} +it \right]^3}
\right\}.
\]
\end{theorem}

\subsection{Higher Codimension Examples.} Here we focus on the case of 
a quadric of codimension two in $\C^4$, given by 
\[
M= \{(z,w) \in \C^2 \times \C^2; \ \Imm { w} = \phi(z,z) \}
\]
where $\phi: \C^2 \times \C^2 \mapsto \C^2 $ is a sesquilinear
form (i.e. $\phi (z, z') = \overline{\phi (z',z)}$).
If the image of the Levi form (i.e. $\phi$) is not contained
in a one-dimensional cone, then  $M$ is biholomorphic to one of the
following three canonical examples (see \cite{Bog91}): 
\begin{itemize}

\item $M_1$ where $\phi(z,z) = (|z_1|^2, |z_2|^2)$

\item $M_2$ where $\phi(z,z) = (2 \Rre (z_1 \bar{z_2}), |z_1|^2 - |z_2|^2 )$

\item $M_3$ where $\phi(z,z) = (2|z_1|^2 , 2  \Rre (z_1 \bar{z_2}) )$

\end{itemize}

The first case, $M_1$, is just the Cartesian product of two Heisenberg groups
in $\C^2$ where solvability in any dimension is not possible.
In the second case, the Levi form of $M_2$
has one positive and one negative eigenvalue
in each totally real direction. Therefore
$\Box_b$ is solvable for $(0,q)$-forms when $q=0, 2$, but not $q=1$
(see \cite{PeRi03}).  
In the third case, the Levi form of $M_3$ has one positive and one negative
eigenvalue for each totally real direction,
except for one exceptional direction where one of the 
eigenvalues is zero. Again, $\Box_b$ is solvable for 
$(0,q)$-forms when $q=0, 2$, but not
when $q=1$ (again see \cite{PeRi03}). In the cases where
solvability is possible, our next theorem provides the
explicit fundamental solution.

\begin{theorem}
\label{highercodim}

\begin{itemize} 

\item For $M_2$, 
the fundamental solution kernel to the Kohn Laplacian $\Box_b$ for $q=0,  \ 2$ is given by 
\[
N(z,t)= \frac{1}{4\pi^3} \frac{1}{(|z|^4 +|t|^2)^{3/2}}.
\]

\item For $M_3$, the fundamental solution kernel to the Kohn Laplacian $\Box_b$ for $q=0$ is given by 
\begin{multline*}
N(z,t) = \frac{1}{\pi^4} \int_0^1 \int_0^{2 \pi} \sigma_1 (\th) \sigma_2 ( \th)
\frac{r^{\sigma_1 -1}}{(1-r^{\sigma_1})(1-r^{\sigma_2)}} \\
\times \frac{2 \, d \th dr}{\left(-i(t_1 \cos \th + t_2 \sin \th)
+E_1 ( \th, r) |z_1|^2 + E_2 ( \th, r) |z_2|^2\right)^3}
\end{multline*}
where
\[
\sigma_1 = \sigma_1 (\th)=1+ \cos \th,  \ \sigma_2 = \sigma_2 (\th) = 1 - \cos \th, \ 
E_j(\th, r) = \frac{\sigma_j (1+r^{\sigma_j})}{1-r^{\sigma_j}}, \ j=1, \ 2.
\]
If $q=2$, then the expression for $N$ is the same except that the factor of 
$r^{\sigma_1 -1}$ is replaced by $r^{\sigma_2 -1}$.

\end{itemize}

\end{theorem}

The integrands in Theorems \ref{mixedeigen} and \ref{highercodim}  
are integrable in $r$ and $r, \ \th$, respectively
when $z_1$ and $z_2$ are not zero, and in each case, 
$N$ can be shown to be locally 
integrable. More detailed estimates of these formulas
in terms of the volumes of nonisotropic balls of the control metric (see \cite{NaStWa85, NaSt06} will be given in a future paper. The outline of this paper
is as follows. Section 3 contains precise definitions of quadrics
and their group structure. Section 4 contains a brief description
of unitary representation theory, which is one of the key tools
which allows us to transform $\Box_b$ to a Hermite operator.
The basic facts about the spectral decomposition of Hermite operators
are given in Section 5. In section 6, we explicitly evaluate 
the spectral decompositions in the cases mentioned in the above theorems.

%
%

\section{Quadric Submanifolds and $\Box_b$} 

%
%
\subsection{Quadric submanifolds}
Let $M$ be the the quadric submanifold in
$\C^n \times \C^m$ defined by
\[
M= \{(z,w) \in \C^n \times \C^m; \ \Imm { w} = \phi(z,z) \}
\]
where $\phi: \C^n \times \C^n \mapsto \C^m $ is a sesquilinear
form (i.e. $\phi (z, z') = \overline{\phi (z',z)}$). For emphasis, we sometimes write $M_\phi$ to denote the dependence
of $M$ on the quadratic function $\phi$. Note that $M_{-\phi}$
is biholomorphic to $M_\phi$ by the change of variables
$(z,w) \mapsto (z,-w)$.

%
%
\subsection{Lie Group Structure.}
By projecting $M \subset \C^n \times \C^m$ onto $G=\C^n \times \R^m$,
the Lie group structure of $M$ is isomorphic to the following
group structure on $G$:
\[
g g'=(z,t)  (z', t') =\big(z+z', t+t' + 2\Imm \phi (z,z') \big).
\]
Note that $(0,0)$ is the identity in this group
structure and that the inverse of $(z,t)$ is $(-z,-t)$.

The \emph{right invariant} vector fields are 
given as follows: let $g \in G$; if $X$ is a vector field,
then we denote its value at $g$ by $X(g)$ as an element
of the tangent space of $M$ at $g$. 
Define
$R_g: G \mapsto G$ by $R_g(g') = g'g$;
a vector field $X$ is {\em right invariant} if and only if
$X(g)=(R_g)_* \{ X(0) \}$, where $(R_g)_* $ denotes
the differential or push forward operator at $g$ as a map
from the tangent space at the origin to the tangent space at $g$. Let 
$v $ be a vector in $\C^n\approx \R^{2n}$
which can be identified with the tangent space
of $M$ at the origin. Let $\partial_v$ be the real
vector field given by the directional 
derivative in the direction of $v$. Then
the right invariant vector field at an arbitrary
$g=(z,w) \in M$ corresponding to $v$ is given by
\[
X_v (g)=  \partial_v + 2 \Imm \phi (v,z) \cdot D_t
=\partial_v - 2 \Imm \phi (z,v) \cdot D_t
\]
(see Section 1 in Peloso/Ricci  \cite{PeRi03}). 
Let $Jv$ be the vector in $\R^{2n}$ which corresponds
to $iv$ in $\C^n$ (where $i=\sqrt{-1}$).
The CR structure on $G$ is then spanned by vectors of the 
form:
\[
Z_v(g)=(1/2)(X_v-iX_{Jv})
=(1/2)(\partial_v - i \partial_{Jv}) -i \bar{\phi (z,v)} \cdot D_t 
\]
and
\[
\bar Z_v(g)=(1/2)(X_v+iX_{Jv})
=(1/2)(\partial_v + i \partial_{Jv}) +i \phi (z,v) \cdot D_t. 
\]
An easy computation gives:
\[
[Z_v, Z_{v'}] = 0 , \ [\bar Z_v, \bar Z_{v'}] =0, \ \ \  
[Z_v, \bar Z_{v'}] = 2i \phi (v, v') \cdot D_t.
\]
The last expression on the right is the Levi form of $M$
as a map from the complex tangent space ($v \in \R^{2n} = \C^n$)
to the totally real directions (spanned by $D_t$).

For 
$\lambda \in \C^m$, let
\[
\phi^\lambda (z,z') = \phi (z,z') \cdot \lambda
\]
where $\cdot$ is the ordinary dot product (without conjugation).
If $\lambda \in \R^m$, then $\phi^\lambda (z,z') $
is a sesquilinear scalar-valued form with an associated Hermitian
matrix. Let $v^\lambda_1, \dots, v^\lambda_n$
be an orthonormal basis for $\C^n$ which diagonalizes this matrix
and we write
\[
\phi^\lambda (v^\lambda_j, v^\lambda_k) = \delta_{jk} \mu_j^\lambda
\]
where $\mulj$, $1 \leq j \leq n$ are its eigenvalues.

\subsection{Special Coordinates.}
\label{special}
For $\lambda\in\R^m$, define the function $\nu(\lambda)$ by
\[
\nu(\lambda) = \text{rank}(\phi^\lambda).
\]
The function $\nu(\lambda)$ satisfies $0\leq \nu(\lambda) \leq n$.
We assume the eigenvalues are ordered so that
$\mu^\lambda_1, \dots , \mu^\lambda_{\nu(\lambda)} \not=0$ and $\mu^\lambda_{\nu(\lambda)+1}, \dots, 
\mu^{\lambda}_n=0$. 
We identify $x$ with $(x^\lambda_1, \dots, x^\lambda_n)$
and $y$ with $(y^{\lambda}_1, \dots, y^{\lambda}_n)$.
We also write $z=\sum_{j=1}^n (x^\lambda_j +i y^\lambda_j) v^\lambda_j$
for $z = x+iy \in \C^n$. 
Additionally, we let $z' = (z_1,\dots, z_{\nu(\lambda)})$, $z'' = (z_{\nu(\lambda)+1},\dots,z_n)$ and 
similarly for $x$ and $y$. Although for many canonical examples, the choice of coordinates will
vary smoothly in $\lambda$, this is not the case in general.

%
%
\subsection{$\Box_b$ Calculations}
Let $v_1, \dots  , v_n$ be any orthonormal basis for $\C^n$
(later, this choice will be a special coordinate basis mentioned just
before Section \ref{special}).
Let $X_j= X_{v_j}$, $Y_j=X_{Jv_j}$,
and let $Z_j=(1/2) (X_j-iY_j)$, $\bar{Z}_j = (1/2) (X_j+iY_j)$
be the right invariant vector fields defined above 
(which are also the left invariant vector fields for the group structure
with $\phi$ replaced by $-\phi$).
Also let $dz_j$ and $d\bar{z}_j$ be the dual basis.
A $(0,q)$-form can be expressed as $\sum_{K\in \I_q} \phi_K\, d\z^K $
where $\I_q = \{ K= (k_1,\dots,k_q) : 1 \leq k_1 < \cdots < k_q \leq n\}$.
We recall the following proposition \cite[Proposition 2.1]{PeRi03}:

\begin{prop}

\begin{equation}
\label{boxb}
\Box_b ( \sum_{K\in\I_q} \phi_K\, d\z^K )
= \sum_{K,L\in \I_q}  \Box_{LK} \phi_K \, d \bar{z}^L
\end{equation}
where
\begin{equation}
\label{boxLK}
\Box_{LK} = - \delta_{LK} \mathcal{L} + M_{LK}
\end{equation}
and $\mathcal{L}$ is the sub-Laplacian on $G$:
\[
\mathcal{L} = (1/2) \sum_{k=1}^n \bar Z_k Z_k + Z_k \bar Z_k
\]
and 
\[
M_{LK} = \left\{
\begin{array}{cc} \displaystyle
\vspace{.05in} \frac 12 \left( \sum_{k \in K} [Z_k , \bar Z_k] - \sum_{k \not\in K} [Z_k , \bar Z_k] 
\right) & \textrm{if} \ K=L \\
\vspace{.1in} \eps(K,L) [Z_k, \bar Z_l] & \textrm{if} \ |K \cap L | = q-1 \\
0 & \textrm{otherwise.}
\end{array}
\right.
\]
Here, $\eps(K,L) $ is $(-1)^{d}$ where 
$d$ is the number of elements in $K \cap L$ between
the unique element $k \in K-L$ and the unique element
$l \in L-K$.

\end{prop}

The above proposition is stated and proved in  \cite{PeRi03} for
left-invariant vector fields. If right invariant
vector fields are used, then the above proposition provides a formula for 
$\Box_b$ associated to $M_{-\phi}$.


For later, we record the diagonal part of $\Box_b$, i.e., 
$\Box_{LL}$. Using (\ref{boxLK}) with $L=K$ and the above formulas for
$Z_k$, we obtain
\begin{equation}
\label{BoxKK}
\Box_{LL} = 
-\frac14 \sum_{k=1}^n \Big(X_k^2 +Y_k^2\Big) 
+i \left( \sum_{k \in L} \phi(v_k, v_k) \cdot D_t 
- \sum_{k \not\in L} \phi(v_k, v_k) \cdot D_t \right).
\end{equation}
For the case of $\Box_b$ on the Heisenberg group,
$\phi(z,z)=|z|^2$ and $Z_k=D_{z_k}-i \bar{z_k} D_t$. In this case,
$\Box_b$ is a diagonal operator (since $[Z_k, \bar{Z}_l] =0$
when $k \not= l$) and the above formula for $\Box_{LL}$
gives the coefficient of $\Box_b$ acting on forms of the type
$\phi_L (z) d \bar{z}^L$.

We will also need the adjoint of $\Box_{LK}$,
which is defined as
\[
\int_{(z,t) \in G} \Box_{LK} \{ f(z,t) \} g(z,t) \, dx\, dy\, dt
=
\int_{(z,t) \in G} f(z,t) \Box^{\textrm{adj}}_{LK} \{g(z,t) \} \, dx\, dy\, dt
\]
(note: this is the ``integration by parts'' adjoint, not 
the $L^2$ adjoint, since there is no conjugation). An easy computation
shows
\begin{equation}
\label{BoxKKadj}
\Box^{\textrm{adj}}_{LL} = 
-\frac14 \sum_{k=1}^n\Big(X_k^2 +Y_k^2\Big) 
-i \left( \sum_{k \in L} \phi(v_k, v_k) \cdot D_t 
- \sum_{k \not\in L} \phi(v_k, v_k) \cdot D_t \right).
\end{equation}
Note the minus sign instead of the plus sign in front of the 
imaginary term at the end.
For later use, note that 
\begin{equation}\label{BoxLLadj and funct form}
\Box_{LL} = \overline{ \Box^\textrm{adj}_{LL} } \ \ 
\textrm{and} \ \ \Box_{LL} \{f(-z,-t) \} = 
( \Box^\textrm{adj}_{LL} f)(-z,-t)
\end{equation}
for any smooth function $f$.

%
%

\section{Representation Theory.}

%
%
\subsection{Unitary Representations}

\begin{definition}

For a Lie group, $G$, a {\em unitary representation} is a 
homomorphism $\pi$ from $G$ to the space of unitary operators
on $L^2(\R^{n'})$ for some $n' \geq 0$.

\end{definition} 

For our quadric 
Lie group $G$, we will fix $\lambda \in \R^m$ as above
and take $n'= \nu(\lambda)$ (the number of nonzero eigenvalues, 
$\mu_j^\lambda$ as in Section \ref{special}). Suppose
$z^\lambda = z=x+iy \in \C^n$ is the special
coordinate system mentioned in Section \ref{special}.
Let $t, \ \lambda \in \R^m$ and $\eta \in \C^{n-\nu(\lambda)}$.
For $g=(x,y,t) \in G$,
define
$\pi_{\lambda,\eta} (x,y,t) :L^2 (\R^{\nu(\lambda)}) \mapsto L^2(\R^{\nu(\lambda)})$
by
\[
\pi_{\lambda,\eta} (x,y,t) (h) (\xi)
= e^{i (\lambda \cdot t + 2\Rre(z''\cdot \bar\eta))} 
e^{-2i \sum_{j=1}^{\nu(\lambda)} \mulj y^\lambda_j ( \xi_j + x^\lambda_j)}
h (\xi +2x')
\]
for $h \in L^2 (\R^{\nu(\lambda)})$ (so $\xi \in \R^{\nu(\lambda)}$). 
Note that if $\eta = \zeta + i\varsigma$, then 
$\Rre(z''\cdot \bar\eta)) = x''\cdot \zeta + y''\cdot \varsigma$.
It is a straight forward 
computation that $\pi$ is a representation for $G$.
For further background information on representation theory (in particular
for the Heisenberg group), see \cite{Tay86}.

If $X$ is a differential operator on $G$ comprised of right invariant
vector fields, then we can compute how $X$ ``transforms'' via 
$\pi_{\lambda,\eta}$ to an operator on $L^2(\R^{\nu(\lambda)})$ 
denoted by $d \pi_{\lambda,\eta} X$, which means
that
\begin{equation}
\label{rightinv}
X_g \{\pi_{\lambda,\eta} (g) \} =d \pi_{\lambda,\eta} X \circ \pi_{\lambda,\eta} (g).
\end{equation}
In words, the $X_g$ on the left side differentiates 
$\pi_{\lambda,\eta} (g)$ with respect to the variable $g$ whereas on the right side,
$d \pi_{\lambda,\eta} X$ is an operator on $L^2(\R^{\nu(\lambda)})$.
The next proposition identifies $d \pi_{\lambda,\eta} (X) $ 
for our basis of right invariant vector fields of $G$.

\begin{prop}
For the right invariant vector fields, $X_j, \ Y_j, \ D_{t_k}$ 
defined in the last section, the following identities hold
as operators on $L^2(\R^{\nu(\lambda)})$:

\begin{align}
\label{righttransfereqn1}
 X_j \{\pi_{\lambda,\eta} \} (g) = d \pi_{\lambda,\eta} X_j
 & = \begin{cases}  2 D_{\xi_j} \circ \pi_{\lambda,\eta} (g)    & 1 \leq j \leq \nu(\lambda) \\
                               2i \zeta_j \circ \pi_{\lambda,\eta} (g) & \nu(\lambda)+1\leq j \leq n      \end{cases} \\
\label{righttransfereqn2}
Y_j \{\pi_{\lambda,\eta} \} (g)= d \pi_{\lambda,\eta} Y_{j}
  &= \begin{cases} -2i \mulj \xi_j  \circ  \pi_{\lambda,\eta} (g)	& 1 \leq j \leq \nu(\lambda) \\
				2i \varsigma_j \circ \pi_{\lambda,\eta} (g) & \nu(\lambda)+1\leq j \leq n      \end{cases} \\
\label{righttransfereqn3}
D_{t_k}\{\pi_{\lambda,\eta} \}(g)  = 
d \pi_{\lambda,\eta} D_{t_k} &= i \lambda_k  \circ  \pi_{\lambda,\eta} (g) \qquad 1\leq k \leq m.
\end{align}

\end{prop}

\n \textbf{Remark.} If we had used left invariant vector fields
instead of right invariant vector fields, then the order of the operators
on the right would have been reversed (i.e. the $D_{\xi_j}$ would appear
on the right of $\pi_{\lambda,\eta} (g) $, etc.). See, for example \cite[equation (8)]{PeRi03}. We prefer the above ordering
of the operators on the right and therefore have chosen to use  right invariant vector fields.

Equation (\ref{righttransfereqn3}) is immediate.
Equations
(\ref{righttransfereqn1}) and (\ref{righttransfereqn2})
are easily shown to hold 
at the origin, $g=0$, since $X_j (0) = D_{x_j}$
and $Y_j (0)= D_{y_j}$; then use right invariance to conclude these
equations hold at all $g \in G$.

Now we compute the ``transform'' of $\Box_{LK}$ and its adjoint,
via $d \pi_{\lambda,\eta}$, using (\ref{BoxKK}) and (\ref{BoxKKadj}).
Note that the coordinates $(z^\lambda_1, \dots, z^\lambda_n)$
were chosen to diagonalize the form $\phi(z,\tilde z) \cdot \lambda$
with eigenvalues $\mu_j^\lambda$.
This observation, together with formulas (\ref{righttransfereqn1}),
(\ref{righttransfereqn2}), and (\ref{righttransfereqn3}) easily 
establish the following proposition

\begin{prop}

\begin{equation}
\label{dbox}
d \pi_{\lambda,\eta} \Box_{LK} = \left\{
\begin{array}{cc}
- \Delta_\xi + |\eta|^2 +  \sum_{j=1}^{\nu(\lambda)} (\mu_j^\lambda)^2  \xi_j^2   -
(\sum_{j\in L} \mu^\lambda_j - \sum_{j \not \in L} \mu^\lambda_j ) &
\textrm{if } \ K=L \\
0 & \textrm{if} \ K \not= L
\end{array}
\right.
\end{equation}

\begin{equation}
\label{dboxadj}
d \pi_{\lambda,\eta} \Box^{\textrm{adj}}_{LK}  = \left\{
\begin{array}{cc}
- \Delta_\xi + |\eta|^2 + \sum_{j=1}^{\nu(\lambda)} (\mu_j^\lambda)^2  \xi_j^2 +(\sum_{j\in L} \mu^\lambda_j - 
\sum_{j \not \in L} \mu^\lambda_j ) &
\textrm{if } \ K=L \\
0 & \textrm{if} \ K \not= L.
\end{array}
\right.
\end{equation}

\end{prop}

Note that $\Box_{LK}$ and its adjoint transform to operators
that only differ by a sign change in the zeroth order terms
involving the $\mu_j^\lambda$.

%
%
\subsection{Group Fourier Transform.}
For $(z,t)\in G$ and fixed $\lambda \in \R^m$, we express $(z,t)=(x,y,t) = (x',y',x'',y'',t) = (x',y', z'',t)$ 
as in Section \ref{special}. The coordinate
$z''$ may be thought of as in $\C^{n-{\nu(\lambda)}}$ or $\R^{2(n-{\nu(\lambda)})}$. 

For an integrable function $f:G \mapsto \C$, we define the 
{\em group Fourier transform} of $f$
as the operator $\pi^{\lambda,\eta}(f): L^2(\R^{\nu(\lambda)}) \mapsto L^2(\R^{\nu(\lambda)})$ where for
$h \in L^2(\R^{\nu(\lambda)})$,
\begin{align*}
\pi^{\lambda,\eta}(f)  \{h\} (\xi)
&=  \int_{(z=x+iy,t) \in G} f(z,t) \pi_{\lambda,\eta} (z,t) (h)(\xi) \, dx\, dy\, dt \\
&=  \int_{(z=x+iy,t) \in G} f(z,t)
e^{i (\lambda \cdot t + 2\Rre(z''\cdot\bar\eta))} e^{-2i \sum_{j=1}^{\nu(\lambda)} \mulj y^\lambda_j ( \xi_j + x^\lambda_j)}
h (\xi +2x')  \, dx\, dy\, dt.
\end{align*}
As before, $x_j$, $y_j$ are the coordinates
for $x,y \in \R^n$ relative to the basis $v^\lambda_1, \dots, v^\lambda_n$.

The following proposition (see \cite[equation (12)]{BoRa11} relates the group transform to the usual Fourier transform
and easily follows from the definition of $\pi^{\lambda,\eta}$.
For us, the usual Fourier transform on $\R^d$ and the notation that we use to denote it is
\[
\hat f(\xi) = f(\hat\xi) = \mathcal F f(\xi) = \frac{1}{(2\pi)^{d/2}} \int_{\R^d} e^{-i x\cdot \xi}f(x)\, dx.
\]

\begin{prop} For $\xi \in \R^{\nu(\lambda)}$ and $h \in L^2 ( \R^{\nu(\lambda)})$,
 
\begin{multline}\label{keyfoureqn}
\pi^{\lambda,\eta}(f)  \{h\}  (\xi)\\ 
=  (2 \pi)^{(2n+m-\nu(\lambda))/2}\int_{x' \in \R^{\nu(\lambda)}}\mathcal{F}_{x'',y,t} \{f (x,y,t) 
e^{-2i \sum_{j=1}^{\nu(\lambda)} \mu^\lambda_j
x_j y_j} \} (x', 2 \mu^\lambda \circ \xi, -2\eta, - \lambda) h(\xi+2x') \, dx'.
\end{multline}
where $\mathcal{F}_{x'',y,t}$ indicates the 
Fourier transform in the $(x'',y,t)$ variables (but not $x'$) and where
$\mu^\lambda \circ \xi = 
(\mu_1 \xi_1 , \dots, \mu_{\nu(\lambda)} \xi_{\nu(\lambda)} )$.

\end{prop}

As with the classical Fourier transform, the above process can
be reversed to identify $f \in L^2(G)$ from its group Fourier transform.
To see this, assume that 
$  \pi^{\lambda,\eta}(f)  $ is known as an operator
on $L^2(\R^{\nu(\lambda)})$ and we wish to identify $f$. For each $a \in \R^{\nu(\lambda)}$,
define $h_a (\xi)= (2\pi)^{-n-m/2} e^{- i\xi \cdot a}$. Set
\[
u^{\lambda,\eta}(a, \xi)  = \pi^{\lambda,\eta}(f)  (h_a)(\xi).
\]
The above definition needs justification since $h_a \not\in L^2(\R^\nu)$.
However, we can multiply $h_a$ by an increasing sequence of cut-off
functions (approaching 1 everywhere) and then take the limit
as our definition of the left side. This technique is carried out carefully in \cite[equation (16)]{BoRa11}. 
Using (\ref{keyfoureqn}), we obtain
\begin{equation}
\label{recover}
u^{\lambda,\eta}(a, \xi) =\pi^{\lambda,\eta}_f (h_a) (\xi) 
=\mathcal{F} \left\{ f (x,y,t) e^{-2i \sum_{j=1}^{\nu(\lambda)} \mu^\lambda_j x_j y_j} \right\} 
(2a, 2\mu^\lambda \circ \xi, - 2 \eta, - \lambda) 
e^{-i \xi \cdot a}
\end{equation}
where $\mathcal{F}$ is the Fourier transform in all variables.  The motivation for the choice of $h 
= h_a$ is that 
it offers the ``missing" exponential needed to relate the full Fourier transform $f$
with $u^{\lambda,\eta}$.
By inverting the Fourier transform
appearing in (\ref{recover}), we obtain the following proposition.

\begin{prop}
\label{findf}
Let $\tilde u^{\lambda,\eta} ( a, b) = u^{-\lambda, -\frac 12\eta} ( a/2, b/(2 \mu^{- \lambda}))$
where $b/(2 \mu^{-\lambda})$ is the vector quantity whose
$j$th component is $b_j/(2 \mu_j^{-\lambda})$. Then
\begin{equation}
\label{inverse}
f (x',y', \w\eta, \w\lambda) = 
e^{-2i \sum_{j=1}^{\nu(\lambda)}  \mulj x_j y_j}
\mathcal{F}^{-1}_{a,b} \left(e^{-\frac i4 \sum_{j=1}^{\nu(\lambda)} a_j b_j /\mulj}
\tilde u^{\lambda,\eta} ( a, b) \right) (x',y')
\end{equation}
where $f (x',y', \w\eta, \w\lambda)$ is the Fourier transform of $f$ in the variables
$\eta$ and $\lambda$  (but not $x'$ and $y'$).
\end{prop}

We can further recover $f$ 
using the inverse Fourier transform in the $\eta, \ \lambda$ variables provided
that the orthonormal basis $v^\lambda_1, \dots, v^\lambda_n$, which diagonalizes
$\phi^\lambda$, depends continuously on $\lambda$ (as it does in the examples
to follow).

\subsection{The Transformed $\Box_b$.}
We start with a general proposition which describes how right invariant 
operators transform via the group transform.

\begin{prop}
\label{propdual}
Suppose $X$ is a differential operator of order at least one which is
comprised of right invariant vector fields. Let $X^{\textrm{adj}}$
denote the ``integration by parts'' adjoint of $X$: i.e., for $f_1\in\Dom(X)\cap L^2(G)$ and $f_2 \in \Dom(X^{\textrm{adj}})\cap L^2(G)$,
\[
\int_{g \in G} Xf_1(g) f_2(g) \,dg =  \int_{g \in G} f_1(g) X^{\textrm{adj}}f_2(g) \,dg
\]
where $dg = dx \, dy \, dt$. Then
\[
d \pi_{\lambda, \eta} (X^{\textrm{adj}}) \circ \pi^{\lambda ,\eta} (f) =  \pi^{\lambda, \eta} (Xf).
\]

\end{prop}

\n \textbf{Proof.} If $h \in L^2 (\R^{\nu(\lambda)})$, then
\begin{align*}
\pi^{\lambda, \eta} (Xf) (h) &= 
\int_{g \in G}Xf(g) (\pi_{\lambda, \eta} (g))(h) \, dg \\
&=\int_{g \in G} f(g) X^{\textrm{adj}}_g (\pi_{\lambda, \eta} (g))(h) \, dg \\
&= \int_{g \in G} f(g) d \pi_{\lambda, \eta} (X^{\textrm{adj}}) \circ
(\pi_{\lambda, \eta} (g))(h) \, dg \\
&= d \pi_{\lambda, \eta}(X^{\textrm{adj}}) 
\Big\{ \int_{g \in G} f(g)\pi_{\lambda, \eta} (g))(h) \, dg \Big\} \\
\end{align*}
which establishes the proposition.

\m

We shall apply this proposition to $\Box_{LL}$
using (\ref{dboxadj}) (note that the off-diagonal
terms transform to zero).
We are interested in finding the 
the operator $N_L$ which satisfies the equation $\Box_{LL} \circ N_L = I-P_L$, where $P_L:L^2(G) \to L^2(G)$
is the orthogonal projection onto the $\Ker(\Box_{LL})$.
For now, we will fix the index $L$.
The operators $N_L$ and $P_L$ are given by group convolution with 
functions, denoted by  $N (z,t) $ and $P (z,t)$, respectively.
So we are interested in solving 
$\Box_{LL} N(z,t) = \delta (z,t) - P(z,t)$ where
$\delta$ is the Dirac delta function supported at the origin.

\begin{corollary}
\label{transfercor}
$\Box_{LL} N(z,t) = \delta (z,t) - P(z,t)$ if and only if
\begin{equation}
\label{transfequiv}
d \pi_{\lambda, \eta} (\Box_{LL}^{\textrm{adj}}) 
(\pi^{\lambda, \eta} N)
= I - P^{\lambda, \eta}
\end{equation}
as operators on $L^2(R^{\nu(\lambda)})$,
where $P^{\lambda, \eta}$ is the $L^2$ projection onto
the kernel of $d \pi^{\lambda, \eta} (\Box_{LL}^{\textrm{adj}})$.

\end{corollary}

\n \textbf{Proof.} We apply Proposition \ref{propdual} to $\Box_{LL}$.
Let $h \in L^2(\R^{\nu(\lambda)})$. The existence of $N$ solving
$\Box_{LL} N(z,t) = \delta (z,t) - P(z,t)$ implies:
\begin{align*}
d \pi_{\lambda, \eta} (\Box_{LL}^{\textrm{adj}}) 
\pi^{\lambda, \eta} N(h) &=
\int_{(z,t) \in G}  \Box_{LL} \{N(z,t) \}  \pi_{\lambda, \eta} (z,t)  h \, dz\, dt \\
&= \int_{(z,t) \in G} (\delta (z,t) - P(z,t)) \pi_{\lambda, \eta} (z,t)  h(\xi) \, dz\, dt \\
&= h(\xi) - \int_{(z,t) \in G} P(z,t) \pi_{\lambda, \eta} (z,t)  h(\xi) \, dz\, dt 
\end{align*}
(since $\pi_{\lambda, \eta}(0,0)$ is the identity operator).
We further claim that the integral on the right is the same
as the $L^2$ projection of $h$ onto the $\Ker\big(d \pi_{\lambda, \eta} (\Box_{LL}^{\textrm{adj}})\big)$.
This is established by proving the following:
\begin{enumerate}
\item $d \pi_{\lambda, \eta} (\Box_{LL}^{\textrm{adj}}) \left\{ \int_{(z,t) \in G} P(z,t) \pi_{\lambda, \eta} (z,t)  h(\xi) \, dz\, dt \right\} =0$
for all $h \in L^2 (G)$

\item if $h \in L^2(G)$ is orthogonal to 
$\Ker\big(d \pi_{\lambda, \eta} (\Box_{LL}^{\textrm{adj}})\big)$, then
\[
\int_{(z,t) \in G} P(z,t) \pi_{\lambda, \eta} (z,t)  h(\xi) \, dz\, dt  =0.
\]
\end{enumerate}

\noindent The first fact is an easy consequence of Proposition \ref{propdual}.
To establish the second fact, it is enough to show
\[
\left \langle \int_{(z,t) \in G} P(z,t) \pi_{\lambda, \eta} (z,t)  h(\xi) \, dz\, dt \ , \ g( \xi) \right\rangle_\xi =0
\]
for all $g $ belonging to the $\Ker\big(d \pi_{\lambda, \eta} (\Box_{LL}^{\textrm{adj}})\big)$.
We have
\[
\left \langle \int_{(z,t) \in G} P(z,t) \pi_{\lambda, \eta} (z,t)  h(\xi) \, dz\, dt \ , \ g( \xi)\right \rangle_\xi
 =\left \langle  h(\xi), \int_{(z,t) \in G} \overline{P(z,t)} \pi^*_{\lambda, \eta} (z,t)  g(\xi) \, dz\, dt\right \rangle_\xi .
\]
Since $h$ is orthogonal to the $\Ker\big(d \pi_{\lambda, \eta} (\Box_{LL}^{\textrm{adj}})\big)$,
it suffices to show
\[
d \pi_{\lambda, \eta} (\Box_{LL}^{\textrm{adj}}) \left\{ \int_{(z,t) \in G} 
\overline{P(z,t)} \pi^*_{\lambda, \eta} (z,t)  g(\xi) \, dz\, dt \right\} =0.
\]
Now $\pi^*_{\lambda, \eta} (z,t) = \pi_{\lambda, \eta} (-z,-t)$, so after a change of variables,
the left side becomes
\begin{multline*}
d \pi_{\lambda, \eta} (\Box_{LL}^{\textrm{adj}}) \left\{ \int_{(z,t) \in G} \overline{P(-z,-t)} \pi_{\lambda, \eta} (z,t)  g(\xi) \, dz dt\right\}\\
= \int_{(z,t) \in G} \overline{P(-z,-t)} (\Box^\textrm{adj}_{LL}   \pi_{\lambda, \eta} (z,t))  g(\xi)  \, dz\, dt \\
= \int_{(z,t) \in G} \overline{P(z,t)} (\Box^\textrm{adj}_{LL}   \pi_{\lambda, \eta}) (-z,-t)  g(\xi)  \, dz\, dt.
\end{multline*}
A chain rule calculation shows that $(\Box^\textrm{adj}_{LL} f)(- z,-t) = \Box_{LL} \{ f(-z,-t) \}$
(recall (\ref{BoxLLadj and funct form})).
Therefore, we can integrate by parts to show that the right side equals
\[
\int_{(z,t) \in G} \Box^\textrm{adj}_{LL}  \{\overline{P(z,t)} \}  \pi_{\lambda, \eta} (-z,-t)  g(\xi) \, dz dt.
\]
Since $ \overline{\Box^\textrm{adj}_{LL}} = \Box_{LL}$, the above term is zero, as desired.
This proves the second fact and
establishes (\ref{transfequiv}). The converse can be established similarly.

\section{The Hermite Operator}

Our starting point is 
equation  Corollary \ref{transfercor}, which allows us to 
transfer the analysis of $\Box_{LL}$ to the following operator
equation
\[
d \pi_{\lambda, \eta} (\Box_{LL}^{\textrm{adj}}) 
\pi^{\lambda, \eta} N
= I - P^{\lambda, \eta}
\]
where according to (\ref{dboxadj}),
\[
d \pi_{\lambda, \eta} (\Box_{LL}^{\textrm{adj}})= - \Delta_\xi + |\eta|^2 + 
\sum_{j=1}^{\nu(\lambda)} (\mu_j^\lambda)^2  \xi_j^2 +(\sum_{j\in L} \mu^\lambda_j - 
\sum_{j \not \in L} \mu^\lambda_j ).
\]
The operator on the right is a modified Hermite operator which has
a well known spectral decomposition which we now describe.
On the real line, and for $\ell$ a nonnegative integer,
define the $\ell^{\textrm{th}}$ Hermite function
\[
\psi_\ell(x) = 
\frac{(-1)^\ell} {2^{\ell/2} \pi^{1/4}(\ell!)^{1/2}}\frac{d^\ell}{dx^\ell}\{ e^{-x^2} \} e^{x^2/2} ,
\ \ x \in \R.
\]
Each $\psi_\ell$ has unit $L^2$-norm on the real line and satisfies
the equation
\[
- \psi_\ell ''(x) + x^2 \psi_\ell (x) = (2\ell+1) \psi_\ell(x),
\]
(see Thangavelu's book \cite{Tha93}). 
The Hermite operator, $-D_{xx} + x^2$ is a nonnegative, self-adjoint operator
and the $\psi_\ell, \ \ell =0, 1, \dots$ forms a complete orthonormal basis
of eigenfunctions with eigenvalues $2\ell+1$.
In several variables, we let $\ell=(\ell_1, \dots , \ell_{\nu(\lambda)})$ where 
each $\ell_i$ is a non-negative integer. For each $\lambda \in R^m \setminus\{0\}$,
define
\[
\psi^\lambda_{\ell_j} (\xi_j) = \psi_{\ell_j} (|\mulj|^{1/2} \xi_j)
|\mulj|^{1/4}; \ \ 
\Psi^\lambda_\ell (\xi) = \prod_{j=1}^{\nu(\lambda)} \psi^\lambda_{\ell_j} ( \xi_j).
\]
Each $\psi^\lambda_{\ell_j} (\cdot )$ has unit $L^2$-norm on $\R$
and hence $\Psi^\lambda_\ell$ has unit $L^2$-norm on $\R^{\nu(\lambda)}$.
An easy calculation shows that 
\begin{equation}
\label{hermite1}
(-D_{\xi_j \xi_j} + (\muj \xi_j)^2) \{ \psi^\lambda_{\ell_j} (\xi_j) \}
= (2\ell_j+1) \psi^\lambda_{\ell_j} (\xi_j) |\muj |.
\end{equation}
Therefore
\begin{align*}
d \pi_{\lambda, \eta} (\Box_{LL}^{\textrm{adj}}) \{\Psi^\lambda_\ell (\xi) \} &=
\Lambda^\lambda_\ell \Psi^\lambda_\ell (\xi) \\
\textrm{where}\ \  \Lambda^{\lambda, \eta}_\ell&= \sum_{j=1}^{\nu(\lambda)} (2 \ell_j+1)|\mu^\lambda_j| +
\left(\sum_{j\in L} \mu^\lambda_j - \sum_{j \not \in L} \mu^\lambda_j \right) + | \eta|^2.
\end{align*}
The collection of functions $\Psi^\lambda_\ell (\xi)$ form a complete set of orthonormal 
eigenfunctions for $d \pi_{\lambda, \eta} (\Box_{LL}^{\textrm{adj}})$ with eigenvalues 
$\Lambda^{\lambda, \eta}_\ell$.
Note that $\Lambda^{\lambda, \eta}_\ell \geq 0$ for all $\ell$ and $\Lambda^{\lambda, \eta}_\ell >0$
for all nonzero multiindices $\ell$. When $\ell=0$,  $\Lambda^{\lambda, \eta}_0$
can equal zero only if $\eta =0$,
all the $\mu^\lambda_k <0$ for $k\in L$, and all $\mu^\lambda_k >0$
for $k \not \in L$. This cannot occur if either the number of negative eigenvalues 
is not equal to $q$=length of ($L$) or if the number of positive
eigenvalues is not equal to $n-q$. 
This condition is the hypothesis of the first part of the following
solvability theorem in Peloso/Ricci \cite{PeRi03}.

\begin{theorem}
\label{existenceth}

\begin{enumerate}

\item (Theorem 1 in \cite{PeRi03}) Suppose $|L|=q$. If the number of negative eigenvalues of $\phi^\lambda$
is not equal to $q$ or if the number of positive
eigenvalues of  $\phi^\lambda$ is not equal to $n-q$, then the 
$\Ker\big(d \pi_{\lambda, \eta} (\Box_{LL}^{\textrm{adj}})\big) = \{0\}$ its inverse is given by the operator
\[
\mathfrak I^{\lambda, \eta} = \sum_\ell \frac{1}{\Lambda^{\lambda, \eta}_\ell} P_\ell^\lambda
\]
where $P_\ell^\lambda$ is the orthogonal projection onto the space
spanned by $\Psi^\lambda_\ell (\xi) $.

\item (Theorem 5.2 in \cite{PeRi03})
Suppose $\nu(\lambda)=n$ (so there is no $\eta$); if $\mu^\lambda_j < 0$  for $j \in L$ and
$\mu^\lambda_j >0$ for $j \not\in L$, then
$\Ker\big(d \pi_{\lambda, \eta} (\Box_{LL}^{\textrm{adj}})\big)$ is the space spanned by 
\[
\Psi^\lambda_0 (\xi) = \prod_{j=1}^n |\mu^\lambda_j|^{1/4} e^{-|\mu^\lambda_j| \xi_j^2/2}.
\]

\end{enumerate}

\end{theorem}

In the first case of this theorem, we will use Proposition \ref{findf}
to identify
$N(z,t)$ (or at least its $t$-Fourier transform) by applying the 
the operator $ \mathfrak I^{\lambda, \eta}$ to the function 
$h_a (\xi) = (2 \pi)^{-n-m/2} e^{-i \xi \cdot a}$, 
as in the discussion leading up to Proposition \ref{findf}.
The first step is to compute
\[
u^{\lambda, \eta} ( a, \xi) = \pi^{\lambda, \eta}(N) \{h_a (\xi) \} = \mathfrak I^{\lambda, \eta} (h_a)( \xi)  
=  (2 \pi)^{-n-m/2}  \sum_{\ell \in \Z^n_+} \frac{1}{\Lambda_\ell^{\lambda,\eta}} \prod_{j=1}^{\nu(\lambda)}  P^\lambda_{\ell_j} \{e^{- i \xi_j  a_j} \}
\]
where $\Z^n_+$ is the set of $n$-tuples of nonnegative integers and
where $P^\lambda_{\ell_j}$ is the orthogonal projection
onto the space of $L^2$ functions in the variable $\xi_j$ spanned
by $\psi^\lambda_{\ell_j} (\xi_j)$.
Each projection term on the right is
\begin{align*}
 P^\lambda_{\ell_j} (e^{- i  \xi_j  a_j} ) 
 &= \left(  \int_{\tilde \xi_j \in \R} e^{-i \tilde \xi_j a_j}
 \psi_{\ell_j} ( |\mu^\lambda_j|^{1/2} \tilde \xi_j)
 |\mu^\lambda_j|^{1/4} \, d \tilde \xi_j \right) |\mu^\lambda_j|^{1/4} 
 \psi_{\ell_j}( |\mu^\lambda_j|^{1/2} \xi_j) \\
 &= (2 \pi)^{1/2}\w{\psi_{\ell_j}} (a_j/|\mu^\lambda_j|^{1/2}) \psi_{\ell_j} ( |\mu^\lambda_j |^{1/2} \xi_j) \\
 &= (2 \pi)^{1/2}(-i)^{\ell_j} 
 \psi_{\ell_j} (a_j/|\mu^\lambda_j|^{1/2}) \psi_{\ell_j} ( |\mu^\lambda_j|^{1/2} \xi_j) 
\end{align*}
where the last equality uses a standard fact that a Hermite function
equals its Fourier transform up to a factor of 
$(-i)^{\ell_j}$ (see \cite{Tha93}). Substituting this expression on the right 
into the above expression for
$u^{\lambda, \eta}( a, \xi)$, we obtain
\begin{equation}
\label{ueqn}
u^{\lambda,\eta} ( a, \xi) = \pi^{\lambda, \eta}(N) \{h_a (\xi) \} =
(2 \pi)^{-n-m/2 +{\nu(\lambda)}/2}  \sum_\ell \frac{(-i)^{\ell}}{\Lambda^{\lambda,\eta}_\ell} 
\prod_{j=1}^{\nu(\lambda)} 
\psi_{\ell_j} (a_j/|\mu^\lambda_j|^{1/2}) \psi_{\ell_j} ( |\mu^\lambda_j|^{1/2} \xi_j).
\end{equation}
In view of (\ref{inverse}), to compute $N (x',y', \w\eta, \w \lambda)$,
we need to determine
\[
\tilde u^{\lambda,\eta} ( a, b) = u^{-\lambda,-\frac 12\eta} ( a/2, b/(2 \mu^{- \lambda}))
\]
where $b/(2 \mu^{-\lambda})$ is the vector quantity whose
$j$th component is $b_j/(2 \mu^{-\lambda}_j)$.
From the previous equality, 
we have
\begin{equation}
\label{utwidd}
\tilde u^{\lambda,\eta} ( a, b) =
(2 \pi)^{-\frac 12(2n+m -\nu(\lambda))}  
\sum_{\ell\in\Z^n_+} \frac{(-i)^{\ell}}{\Lambda^{-\lambda,-\frac 12\eta}_\ell}
\prod_{j=1}^{\nu(\lambda)}
\psi_{\ell_j} (a_j/2|\mulj|^{1/2}) \psi_{\ell_j} (b_j |\mulj|^{1/2}/2 \mu^{-\lambda}_j).
\end{equation}

Using (\ref{inverse}),
the formula for the partial transform of $N$ is given in the 
following proposition.

\begin{prop} 

The partial $(z'',t)$-Fourier transform of the 
fundamental solution to $\Box_L$ is given by
\begin{equation}
\label{Nform}
N (x',y', \w\eta, \w\lambda) = 
e^{-2i \sum_{j=1}^{\nu(\lambda)}  \mulj x_j y_j}
\mathcal{F}^{-1}_{a,b} \left(e^{-\frac i4 \sum_{j=1}^{\nu(\lambda)} a_j b_j /\mulj}
\tilde u^{\lambda,\eta} ( a, b) \right) (x',y')
\end{equation}
where $\tilde u^{\lambda,\eta} ( a, b)$ is given in 
(\ref{utwidd}).

\end{prop}

In subsequent sections, this formula will be explicitly computed
in the examples mentioned in the introduction.

In a similar fashion, the Szeg\"o operator, $P$, representing the projection
onto the zero eigenspace in the second case of Theorem \ref{existenceth}
can be computed. We obtain the following result.

\begin{prop}

Suppose $\nu(\lambda)=n$ (so there is no $\eta$); 
if $\mu^\lambda_j < 0$  for $j \in L$ and
$\mu^\lambda_j >0$ for $j \not\in L$, then the operator representing
the projection onto the zero-eigenspace of $\Box_L$ is given 
by a (group) convolution with the kernel $P(z,t)$ whose $t$-Fourier
transform is given by
\[
P(x,y,\hat \lambda) = (2 \pi)^{-(n+m/2)}   \prod_{j=1}^n
|\mu^\lambda_j| e^{-|\mu^\lambda_j| (x_j^2+y_j^2)}.
\]

\end{prop}

This can be easily inverted using the inverse Fourier transform
in $\lambda$ to recover the classical formulas for the Szego kernel.


\section{Calculation of N} 

\subsection{Reductions for the General Case.}
Our goal is now to
try to unravel the formula (\ref{Nform})
for $N(x,y, \w\eta, \w\lambda)$ under the hypothesis
in the first
part of Theorem \ref{existenceth}. 
The discussion in this section will apply to the general
case. Subsequent sections will address the specific cases given in 
Theorems \ref{mixedeigen} and \ref{highercodim}.

We rewrite (\ref{Nform}) as
\begin{align*}
N(z', \w\eta, \w\lambda) &= (2 \pi)^{-\frac 12(2n+m -\nu(\lambda))} e^{-2i \sum_{j=1}^n \mu^\lambda_j x_j y_j } \\
& \sum_{\ell \in Z^n_+} \frac{(-i)^\ell}{\Lambda^{-\lambda,-\frac \eta2}_\ell}
\prod_{j=1}^{\nu(\lambda)}  \mathcal{F}^{-1}_{a_j} \left\{
\psi_{\ell_j} (\frac{a_j}{2  |\mu^\lambda_j|^{1/2}})
\mathcal{F}^{-1}_{b_j} \left(e^{(-i/4)  \frac{a_j b_j}{\mu^\lambda_j}}
\psi_{\ell_j}( \frac{b_j |\mu^\lambda_j|^{1/2}}{-2 \mu^\lambda_j })\right)(y_j)
\right\}(x_j) .
\end{align*}
Here, we are using the following notation for the one-variable
Fourier transform:
\[
\mathcal{F}^{-1}_b (g)(y) = 
\frac{1}{(2 \pi)^{1/2}}\int_{-\infty}^\infty g(b)e^{i y b} \, db
\]
and similarly for $\mathcal{F}^{-1}_a (g)(x)$.
We compute the above partial inverse Fourier transform expression in $a_j$:
\begin{align*}
&\mathcal F_a^{-1}\bigg\{ \psi_{\ell_j}\Big(\frac{a}{2|\mu|^{1/2}}\Big) \mathcal F_b^{-1} \bigg( e^{-i \frac{ab}{4\mu}} \psi_{\ell_j}\Big(\frac{b|\mu|^{1/2}}{-2\mu}\Big) \bigg)(y)\bigg\}(x) \\
&= 4|\mu| \mathcal F_a^{-1}\bigg\{ \psi_{\ell_j}(a)  \psi_{\ell_j}\Big(a - 2\frac{\mu}{|\mu|^{1/2}}y\Big)\bigg\}(2|\mu|^{1/2} x) \\
&=4|\mu|(-1)^{\ell_j} \mathcal F^{-1}\bigg\{ \mathcal F \psi_{\ell_j}(a) \mathcal F \Big\{ e^{i 2\frac{\mu}{|\mu|^{1/2}} y a }\psi_{\ell_j}(a)\Big\}\bigg\}(2|\mu|^{1/2} x)\ \\
&= \frac{4|\mu|}{\sqrt{2\pi}} (-1)^{\ell_j} \int_\R \psi_{\ell_j}(a) e^{i2\frac{\mu}{|\mu|^{1/2}} y(2|\mu|^{1/2}x-a)}\psi_{\ell_j}(2|\mu|^{1/2} x-a)\, da \\
&= \frac{4|\mu|}{\sqrt{2\pi}}  (-1)^{\ell_j} e^{i4\mu xy} \int_\R  \psi_{\ell_j}(a) e^{-2ia y\mu/|\mu|^{1/2} } \psi_{\ell_j}(2|\mu|^{1/2}x-a)\, da.
\end{align*}
Thus, we have the following lemma.
\begin{lemma} 
\label{lemmaequal}
In the first case of Theorem \ref{existenceth}
where the eigenvalues, $\Lambda^{\lambda, \eta}_\ell$, are nonzero
for $\lambda \not=0$,
\label{genformula}
\begin{multline*}
N(z, \w\eta, \w\lambda) =
(2\pi)^{-n-\frac m2} 4^{\nu(\lambda)} e^{2i \sum_{j=1}^n \mu^\lambda_j x_jy_j}  \\
\sum_{\ell} \frac{(-1)^{|\ell|}}{\Lambda^{-\lambda,-\frac \eta2}_\ell}
\prod_{j=1}^{\nu(\lambda)} \int_{\R} |\mu^\lambda_j|\psi_{\ell_j}(a_j)
\psi_{\ell_j} (2 |\mu^\lambda_j|^{1/2} x_j-a_j)
e^{-2i \mu^\lambda_j y_j a_j/|\mu^\lambda_j|^{1/2}} \, d a_j \\
\end{multline*}
where 
\[
\Lambda^{-\lambda,-\frac \eta2}_\ell = \sum_{j=1}^n (2 \ell_j+1) |\mu^\lambda_j| -
\left( \sum_{k \in L} \mu^\lambda_k - \sum_{k \not\in L} \mu^\lambda_k \right) 
+ \frac{|\eta|^2}4.
\]
\end{lemma}

\subsection{The Case when the Eigenvalues have the same Absolute Value.}

In this section, we assume $\nu(\lambda)=n$ (so there is no $\eta$ variables)
and that all the $\mu^\lambda_j$ have the same absolute value, for each $\lambda \in \R^m$.
In this case, we can rescale and assume
\[
|\mu^\lambda_j| = |\lambda|, \ \ \textrm{for all} \  1 \leq j \leq n
\]
and write 
\[
\mu^\lambda_j = \sigma^\lambda_j |\lambda| \ \ \textrm{where} \ 
\sigma^\lambda_j=\pm 1.
\]
This assumption, though restrictive, will allow us
to recover $N$ for the Heisenberg group, 
as well as for the case of $M_2$ with real codimension two,
given in Theorem \ref{highercodim}. 
We will further assume there are no zero eigenvalues
as in the first case of Theorem \ref{existenceth}.
In this case, we have
\[
\Lambda^{-\lambda}_\ell = 2(\sum_{j=1}^n  \ell_j +J)|\lambda|
\]
where $J =J(\lambda)$ is a positive integer.
In the case of the Heisenberg group with $\phi(z,z)=-|z|^2$ and $L$ is an index of length $q$, 
then $J=q$ if $\lambda <0$ and $J=n-q$ if $\lambda >0$.
If $1 \leq q \leq n-1$, then $J$ is a positive integer. 

From Lemma \ref{lemmaequal}, $N$ can be written as 
\begin{equation}
\label{N1}
N(z, \hat \lambda) = 
(2\pi)^{-n-\frac m2}4^n |\lambda|^{n-1} e^{2i (\sigma x \cdot y) |\lambda|}
\sum_{k=0}^\infty \frac{(-1)^k}{2(k+J)}
\sum_{|\ell|=k} \int_{a \in \R^n}\Psi_\ell (a) \Psi_\ell (2 |\lambda|^{1/2} x -a) e^{-2i (\sigma y \cdot a)
|\lambda|^{1/2}}
\end{equation}
where
\[
\Psi_\ell(t) = \psi_{\ell_1} (t_1) \cdots \psi_{\ell_n} (t_n)
\ \ \textrm{and} \ \sigma x = (\sigma^\lambda_1  x_1, \dots, \sigma^\lambda_n
x_n).
\]

Now we use Mehler's formula (see \cite{Tha93}), for each fixed $x, \ y \in \R^n$
\[
\sum_{k=0}^\infty (-r)^k \sum_{|\ell| =k} \Psi_\ell (x) \Psi_\ell (y)
= \phi(-r) , \ \ |r|<1
\]
where 
\[
\phi(r) = \frac{1}{\pi^{n/2} (1-r^2)^{n/2}} e^{-(\frac{1-r}{1+r}) (x+y)^2/4 -
(\frac{1+r}{1-r}) (x-y)^2/4}.
\]
The key idea is to now multiply 
Mehler's formula by $r^{J-1}$ and integrate $r$ over the interval $0 \leq r <1$:
\[
\sum_{k=0}^\infty \frac{(-1)^k}{k+J} \sum_{|\ell| =k} \Psi_\ell (x) \Psi_\ell (y)
= \int_0^{1} r^{J-1} \phi(-r)\, dr.
\]
Note that there is no problem with convergence of the integral
on the right since $J$ is a positive integer.
Applying this formula to (\ref{N1}) with $x$ replaced by $a$
and $y$ replaced by $2 |\lambda|^{1/2} x -a$ gives
\begin{multline*}
N(z, \hat \lambda) = \frac{(2\pi)^{-n-\frac m2} 4^n |\lambda|^{n-1}}{2 \pi^{n/2}}
e^{2i (\sigma x \cdot y) |\lambda|} \\
  \int_{a \in \R^n} \int_0^{1}
\frac{r^{J-1}}{(1-r^2)^{n/2}} e^{-2i(\sigma y \cdot a) |\lambda|^{1/2}}
e^{-(\frac{1+r}{1-r}) |\lambda| x^2} e^{-(\frac{1-r}{1+r}) (|\lambda|^{1/2} x -a)^2}
\, da\, dr.
\end{multline*}
Now integrate out $a$ (completing the square in the exponential, etc.)
to obtain
\[
N(z, \hat \lambda) = \frac{(2\pi)^{-n-\frac m2} 4^n |\lambda|^{n-1}}{2}
\int_0^{1} \frac{ r^{J-1}}{(1-r^2)^{n/2}}
\left(\frac{1+r}{1-r}\right)^{n/2} e^{-(\frac{1+r}{1-r}) |z|^2 |\lambda|}
\, dr.
\]
If $z \not=0$, then the above integral converges.
Now 
set $s=(1+r)/(1-r)$ and then translate $s$ by one unit and we obtain the
following lemma.

\begin{lemma} In the case where $J>0$, as above
\label{equalmu}
\[
N(z, \hat \lambda) = (2\pi)^{-n-\frac m2} 2^n |\lambda|^{n-1}
\int_0^\infty s^{J-1} (s+2)^{n-J-1} e^{-(s+1) |\lambda| |z|^2}
\, ds.
\]
\end{lemma}

This integral can be computed using integration by parts.
$N$ can then be computed using the inverse Fourier transform 
in the $\lambda$ variable, provided the coordinates $z^\lambda$
from Section \ref{special} vary continuously in $\lambda$.
In the Heisenberg group (where $\phi(z,z) = -|z|^2)$, the classical 
formulas for $N$ (see \cite{Ste93}) can be determined by computing the inverse
Fourier transform in $\lambda$ and separating out the integral over $\lambda >0$,
where $J=q$ and $\lambda <0$, where $J=n-q$.

\subsection{Proof of Theorem \ref{highercodim} for $M_2$}

For $\lambda \in \R^2$, 
let $\phi^\lambda: \C^2 \times \C^2 \to \R$ be defined by
$\phi^\lambda (z, z') = \phi (z, z') \cdot \lambda$ 
where $\phi$ is the defining function for $M_2$
given in Theorem \ref{highercodim}. It is easy to compute
that the Hermitian form associated with $\phi^\lambda$ has
two eigenvalues: $\mu^\lambda_1 = |\lambda|$ 
and $\mu^\lambda_2 = - |\lambda|$. 
For $\lambda=|\lambda|( \cos \th, \sin \th ) \in \R^2$,
the associated eigenvectors are
$v^\lambda_1 = \frac 1{ \sqrt{2}} \left(
\frac{\cos  \theta}{\sqrt{1 - \sin \theta} }, \ \sqrt{1- \sin \theta }
\right)$ and $v^\lambda_2 = \frac{1}{\sqrt 2} \left(
-\sqrt{1- \sin \theta }, \ \frac{\cos\theta}{\sqrt{1-\sin\theta}} \right)$, 
respectively.
These eigenvectors vary smoothly with $\th$.

Since the eigenvalues have the same absolute value, we can compute 
$N$ using the techniques in the previous section.
Solvability of $\Box_b$ 
is expected for $(0,q)$-forms where $q=0,2$, but not for $q=1$.
For both $q=0,2$, $\Lambda^{-\lambda}= 2(\ell_1+\ell_2 +1)$
and so we can use Lemma \ref{equalmu} with $J=1$ to 
compute
\[
N(z, \hat \lambda) = \frac{1}{2\pi^3} |\lambda|
\int_0^\infty e^{-(s+1)|\lambda||z|^2} \, ds
= \frac{2 e^{-|\lambda||z|^2}}{\pi^2|z|^2}.
\]
Using the inverse Fourier transform in $t \in \R^2$, we obtain
\[
N(z,t) = \frac{1}{4\pi^4|z|^2}
\int_{\lambda \in \R^2}  e^{-|\lambda||z|^2 +i \lambda \cdot t}
\, d \lambda.
\]
Polar coordinates can be used to reduce the above integral to
\[
N(z,t) = \frac{1}{ 4\pi^4|z|^2} \int_0^{2 \pi}
\frac{d \theta}{(|z|^2-i(t_1 \cos \theta +t_2 \sin \theta))^2}.
\]
The above integrand is periodic in $\theta$. A shift in $\theta$ (specifically, $\theta \mapsto \theta - A$ where $\cos A = t_1/|t|$ and $\sin A = t_2/|t|$) can 
be used to reduce the denominator of the integrand to 
$(|z|^2-i(|t| \cos \theta))^2 $. From here, a residue calculation 
(or Maple) gives
\[
N(z,t) = \frac{1}{2\pi^3}\frac{1}{ (|z|^4+|t|^2)^{3/2}},
\]
as stated in Theorem \ref{highercodim}.
This concludes the proof of the first part of Theorem \ref{highercodim}.
The second part (for $M_3$) will be given after the next section
where we introduce techniques for handling eigenvalues which are not
equal in absolute value.

\subsection{Proof of Theorem \ref{mixedeigen}}
Here, $n=2$ and the multiindex $L=(1,0)$. The eigenvalues are
$\mu_1^\lambda = \sigma_1 \lambda$,
$\mu_2^\lambda = \sigma_2 \lambda$ with $\sigma_1, \sigma_2>0$.
We first consider the case when $\lambda>0$ and we obtain
\[
\Lambda^{-\lambda}_{\ell_1, \ell_2} = 2(\ell_1 |\mu_1|+(\ell_2 +1)|\mu_2|)
= 2 \lambda (\ell_1 \sigma_1+(\ell_2 +1)\sigma_2)
\]
(in this case, there is some cancellation in $\mu_1$-terms
in the formula for $\Lambda^{-\lambda}_{\ell}$ but not 
in $\mu_2$).
From Lemma \ref{genformula}, the operator $N$ becomes 
\begin{equation}
\label{N2}
N(z, \hat \lambda) = 
\frac{8}{(2\pi)^{5/2}} \lambda e^{2i (\sigma x \cdot y) \lambda}
\sum_{\ell_1, \ell_2} \frac{(-1)^{\ell_1+\ell_2}
\sigma_1 \sigma_2}{(\ell_1 \sigma_1+(\ell_2 +1)\sigma_2)}
\int_{a \in \R^2} E(a, x, y, \lambda) \, da
\end{equation}
where
\[
E(a, x, y, \lambda) = 
\prod_{j=1}^2 \psi_{\ell_j} (a_j) \psi_{\ell_j} 
(2 \sigma_j^{1/2}\lambda^{1/2} x_j -a_j) e^{-2i 
\sigma_j^{1/2} y_j  a_j 
\lambda^{1/2}}.
\]
This time, we use Mehler's formula with fractional powers of $r$:
\begin{multline*}
 \sum_{\ell_1, \ell_2} (-r^{\sigma_1})^{\ell_1}(-r^{\sigma_2})^{\ell_2} \psi_{\ell_1} (X_1) 
\psi_{\ell_1} (Y_1)\psi_{\ell_2} (X_2)\psi_{\ell_2} (Y_2) \\
= \frac{1}{\pi}\prod_{j=1}^2 \frac{1}{(1-r^{2\sigma_j})^{1/2}}  
e^{-\left(\frac{1+r^{\sigma_j}}{1-r^{\sigma_j}}\right)
(X_j+Y_j)^2/4 - \left(\frac{1-r^{\sigma_j}}{1+r^{\sigma_j}}\right)
(X_j-Y_j)^2/4} \\
\end{multline*}
where $X_j, \ Y_j \in \R$.
Multiplying this expression by $r^{\sigma_2 -1}$, integrating from $r=0$ to $r=1$ 
and setting $X_j=a_j$ and $Y_j=2\lambda^{1/2}\sigma_j^{1/2} x_j -a_j$, $j=1,2$, gives
\begin{multline*}
 \sum_{\ell_1, \ell_2} \frac{(-1)^{\ell_1+\ell_2}}{\ell_1 \sigma_1+(\ell_2 +1)\sigma_2}
\prod_{j=1}^2
\psi_{\ell_j} (a_j) \psi_{\ell_j} (2\lambda^{1/2} \sigma_j^{1/2} x_j -a_j) 
e^{-2i (\sigma_j^{1/2}  y_j  a_j) \lambda^{1/2}}\\
= \frac 1\pi\int_0^1 \frac{r^{\sigma_2-1}}{(1-r^{2\sigma_1})^{1/2} (1-r^{2\sigma_2})^{1/2}}
\prod_{j=1}^2 e^{-\lambda\sigma_j \frac{1+r^{\sigma_j}}{1-r^{\sigma_j}}x_j^2
- \frac{1-r^{\sigma_j}}{1+r^{\sigma_j}} (\lambda^{1/2}\sigma^{1/2}_j x_j-a_j)^2 -
2i\sigma_j^{1/2} y_j a_j \lambda^{1/2}} \, dr.
\end{multline*}
We then integrate out $a \in \R^2$ (completing the squares in the exponent, etc.) and simplify to 
obtain
\begin{equation}
\label{Nnearfinal}
N(z, \hat \lambda) =
\frac{8}{(2 \pi)^{5/2}} \lambda \sigma_1\sigma_2\int_0^1 \frac{r^{\sigma_2-1}}{(1-r^{\sigma_1}) (1-r^{\sigma_2})}
\prod_{j=1}^2 e^{-\lambda \sigma_j| \frac{1+r^{\sigma_j}}{1-r^{\sigma_j}} |z_j|^2} \, dr.
\end{equation}

When $\lambda$ reverses sign, then so do the $\mu^{\lambda}_j$. 
As a result, the $\mu^\lambda_2$ cancels instead of the $\mu^\lambda_1$
in the expression
for $\Lambda^{-\lambda}_\ell$. 
The above computation goes through unchanged except that $r^{\sigma_1 -1}$
replaces $r^{\sigma_2-1}$  and $|\lambda|$ replaces $\lambda$. We then obtain
\[
N(z, \hat \lambda) =
\frac{8}{(2 \pi)^{5/2}}|\lambda| \sigma_1\sigma_2 \int_0^1 \frac{r^{\sigma_1-1}}{(1-r^{\sigma_1}) (1-r^{\sigma_2})}
\prod_{j=1}^2 e^{-|\lambda| \sigma_j \frac{1+r^{\sigma_j}}{1-r^{\sigma_j}} |z_j|^2} \, dr.
\]

We can evaluate the inverse Fourier transform in $\lambda$ by computing two integrals: one
where $\lambda >0$ and the other where $\lambda<0$.
After a simple integration by parts we obtain:
\begin{multline*}
N(z,t)= \frac{1}{\pi^3} \int_0^1 \frac{\sigma_1\sigma_2}{(it+s_1(r)\sigma_1|z_1|^2+
s_2(r)\sigma_2|z_2|^2)^2}
\frac{r^{\sigma_1-1}\, dr}{(1-r^{\sigma_1})(1-r^{\sigma_2})} \\
+ \frac{1}{\pi^3} \int_0^1 \frac{\sigma_1\sigma_2}{(-it+s_1(r)\sigma_1|z_1|^2
+s_2(r)\sigma_2|z_2|^2)^2}
\frac{r^{\sigma_2-1}\, dr}{(1-r^{\sigma_1})(1-r^{\sigma_2})}
\end{multline*}
where 
\[
s_j(r) = \frac{1+r^{\sigma_j}}{1-r^{\sigma_j}} \ \ \textrm{for} \ j=1, \ 2.
\]
This completes the proof of Theorem \ref{mixedeigen}. 
\m

\n \textbf{Remark.} The same process can be used to evaluate
$N(z,t)$ in cases where $n>2$ or where the eigenvalues
of the Levi form (the  $\sigma_j$) have both positive
and negative terms. 

\subsection{Proof of Theorem \ref{zeroeigen}}
The quadric hypersurace of interest in Theorem \ref{zeroeigen} is 
\[
M= \{ (z_1, z_2, z_3, w) \in \C^4; \ \Rre w = |z_1|^2 +|z_2|^2 \}.
\]
Note that this quadric has a Levi form with diagonal entries $1, \ 1, \ 0$
corresponding to the directions $z_1, \ z_2, \ z_3$, respectively.
The $(0,1)$ forms under consideration are spanned by 
$d \bar z_1$ and $d \bar z_2$ (but not $d \bar z_3$).
Our starting point is the formula for $N (x',y', \w\eta, \w\lambda)$
given in  Lemma \ref{lemmaequal} with $n=3, \ m=1$, $\nu(\lambda) = 2$.
In addition,
$z' = x'+iy'=(z_1, z_2)$; $\eta$ is the Fourier transform
variable associated to the zero-eigendirection variable $z_3$;
$\lambda$ is the Fourier transform variable associated
with $t= \Rre w$; and $\mu^\lambda_1=\mu^\lambda_2 = \lambda$.
We therefore have
\[
\Lambda^{-\lambda,-\eta/2}_\ell = 2|\lambda| (|\ell| +1) + \frac{|\eta|^2}{4}
\]
(note that the term 
$\sum_{j \in L} \mu^\lambda_j -\sum_{j \not \in L} \mu^\lambda_j$ 
is zero since $\mu^\lambda_j = \lambda$, $j=1,2$ and one of 
these indices belongs to $L$ and one does not).

By Lemma \ref{lemmaequal},
\begin{multline*}
N(z', \hat \eta, \hat \lambda)= 
\frac{16|\lambda|^2}{(2 \pi)^{7/2}}  e^{2i \sum_{j=1}^2 \lambda x_j y_j} \\
\sum_{k=0}^\infty \frac{(-1)^k}{2|\lambda|(k+1) + |\eta|^2/4}
\sum_{|\ell|=k} \prod_{j=1}^2 
\int_{\R} \psi_{\ell_j}(a_j)
\psi_{\ell_j} (2 |\lambda|^{1/2} x_j-a_j)
e^{-2i \lambda y_j a_j/|\lambda|^{1/2}} \, d a_j .\\
\end{multline*}

Now the idea is to write
\[
\frac{(-1)^k}{2|\lambda|(k+1) + |\eta|^2/4}
= \frac{(-1)^k}{2|\lambda| \left(k+1+\frac{|\eta|^2}{8|\lambda|} \right)}
=\frac{1}{2|\lambda|} \int_0^1 (-1)^k r^{k+\frac{|\eta|^2}{8|\lambda|}} \, dr.
\]
The computations in the previous section  (using Mehler's formula
and integrating out $a$) can then be repeated with
minor modifications to show the following formula for $N(z', \hat \eta, \hat \lambda)$
which is analogous to (\ref{Nnearfinal}):
\[
N(z', \hat \eta, \hat \lambda)= \frac{8|\lambda|}{(2 \pi)^{7/2}}
\int_0^1 e^{- |\ln r| \frac{|\eta|^2}{8|\lambda|}}
\frac{1}{(1-r)^2} e^{- |\lambda| |z'|^2 \left( \frac{1+r}{1-r} \right)} \, dr.
\]

Note that the above expression is a Gaussian in $\eta$ and so its inverse
Fourier transform in $\eta$ is easily calculated:
\[
N(z', z_3, \hat \lambda) =  \frac{32|\lambda|^2}{(2 \pi)^{7/2}}
 \int_0^1 \frac{1}{|\ln r| (1-r)^2}
e^{- |\lambda| \left( |z'|^2 \left( \frac{1+r}{1-r} \right) + |z_3|^2 \frac{2}{|\ln r|}
\right)}
\, dr.
\]

The inverse Fourier transform in $\lambda$ is now easily calculated:
\[
N(z,t) = \frac{8}{\pi^4}
\int_0^1 \frac{dr}{ |\ln r| (1-r)^2} \textrm{Re} \left\{
\frac{1}{\left[ |z'|^2 \left(\frac{1+r}{1-r}\right) + |z_3|^2 \frac{2}{ |\ln r|} +it \right]^3}
\right\}
\]
which establishes Theorem \ref{zeroeigen}.

\subsection{Proof of Theorem \ref{highercodim} for $M_3$}
For $\lambda \in \R^2$, 
let $\phi^\lambda: \C^2 \times \C^2 \to \R$ be defined by
$\phi^\lambda (z, z') = \phi (z, z') \cdot \lambda$ 
where $\phi$ is the defining function for $M_3$
given in Theorem \ref{highercodim}. We write $\lambda \in \R^2$
in polar coordinates $\lambda= |\lambda| (\cos \th, \sin \th) $.
It is easy to compute
that the Hermitian form $\phi^\lambda$ has
two eigenvalues: $\mu^\lambda_1 = |\lambda|(\cos \th +1)$ 
and $\mu^\lambda_2 =  |\lambda| ( \cos \th -1)$. 
The associated eigenvectors are 
$v^\lambda_1 = \frac{1}{\sqrt 2}
\left(\frac{\sin\theta}{\sqrt{1-\cos\theta}} , \ \sqrt{1 - \cos \theta} \right)$,
and 
$v^\lambda_2 = \frac{1}{\sqrt 2}
\left( -\sqrt{1 - \cos \theta} , \frac{\sin\theta}{\sqrt{1-\cos\theta}}
\right)$, which vary smoothly in $\th$.
Since the 
eigenvalues have opposite sign for all $\th$ except $\th =0$,
we expect solvability of $\Box_b$ on $(0, q)$-forms for $q=0,2$.
We first assume $q=0$. We see that 
\[
\Lambda^{-\lambda}_{\ell_1, \ell_2} = (2 \ell_1+1)| \mu^\lambda_1|
+(2 \ell_2 +1) | \mu^\lambda_2| + \mu^\lambda_1 + \mu^\lambda_2 = 2|\lambda| \left((\ell_1+1) \sigma_1 + \ell_2 \sigma_2 \right)
\]
where $\sigma_1 = (1+ \cos \th)$ and $\sigma_2 = (1- \cos \th)$.
Repeating the same arguments from the last section,
we obtain
\[
N(z, \w \lambda) = \frac{\sigma_1 \sigma_2 |\lambda|}{\pi^3}
\int_0^1 \frac{r^{\sigma_1-1}}{(1-r^{\sigma_1}) (1-r^{\sigma_2})}
\prod_{j=1}^2 e^{-|\lambda| \sigma_j \frac{1+r^{\sigma_j}}{1-r^{\sigma_j}} |z_j|^2} \, dr.
\]
We now take the inverse Fourier transform of this expression in $\lambda \in \R^2$
using polar coordinates. Integrating $|\lambda|$ from $0$ to infinity is a straight
forward integration by parts. However, the $\theta$ and $r$ integrals cannot be
evaluated in closed form. The result is 
\begin{multline*}
N(z,t) = \frac{1}{\pi^4} \int_0^1 \int_0^{2 \pi} \sigma_1 (\th) \sigma_2 ( \th)
\frac{r^{\sigma_1 -1}}{(1-r^{\sigma_1})(1-r^{\sigma_2)}} \\
\times \frac{2 \, d \th dr}{\left(-i(t_1 \cos \th + t_2 \sin \th)
+E_1 ( \th, r) |z_1|^2 + E_2 ( \th, r) |z_2|^2\right)^3}
\end{multline*}
where
\[
\sigma_1 = \sigma_1 (\th)=1+ \cos \th,  \ \sigma_2 = \sigma_2 (\th) = 1 - \cos \th, \ 
E_j(\th, r) = \frac{\sigma_j (1+r^{\sigma_j})}{1-r^{\sigma_j}}, \ j=1, \ 2
\]
as stated in Theorem \ref{highercodim}. In the case where $q=2$, 
$\Lambda^{-\lambda}_{\ell_1, \ell_2}$ becomes
$2|\lambda| \left(\ell_1 \sigma_1 +(\ell_2+1) \sigma_2 \right)$.
This change results
in a factor of $r^{\sigma_2 -1}$ in the numerator of $N$ instead of the factor
of $r^{\sigma_1 -1}$. This completes the proof of Theorem \ref{highercodim}.

\bibliographystyle{alpha}
\bibliography{mybib}

\end{document}